# On Statistics of the Riemann Zeros Differences

Yuri Bachilov*

*Analogic Corporation, CT Medical Division, Peabody, MA 01960*

Abstract

Numerical study of the distribution of the Riemann zeros differences following the work [1] shows the significance of the function for which the prime sum expression is proposed. Computational results related to this definition explored with various prime cut-offs.

* ybachilov@gmail.com

## Approach

This paper is devoted to the investigation of the remarkable observation of Ricardo Pérez Marco [1] named as "Riemann zeros repel their deltas".

Fig. 1 represents section of two histograms (appropriately scaled, obtained with the spacing =0.001) of RZ deltas for the first $5{\cdot}10^6$(5M) and $10^9$ (1B) of Riemann zeros[1] [6].

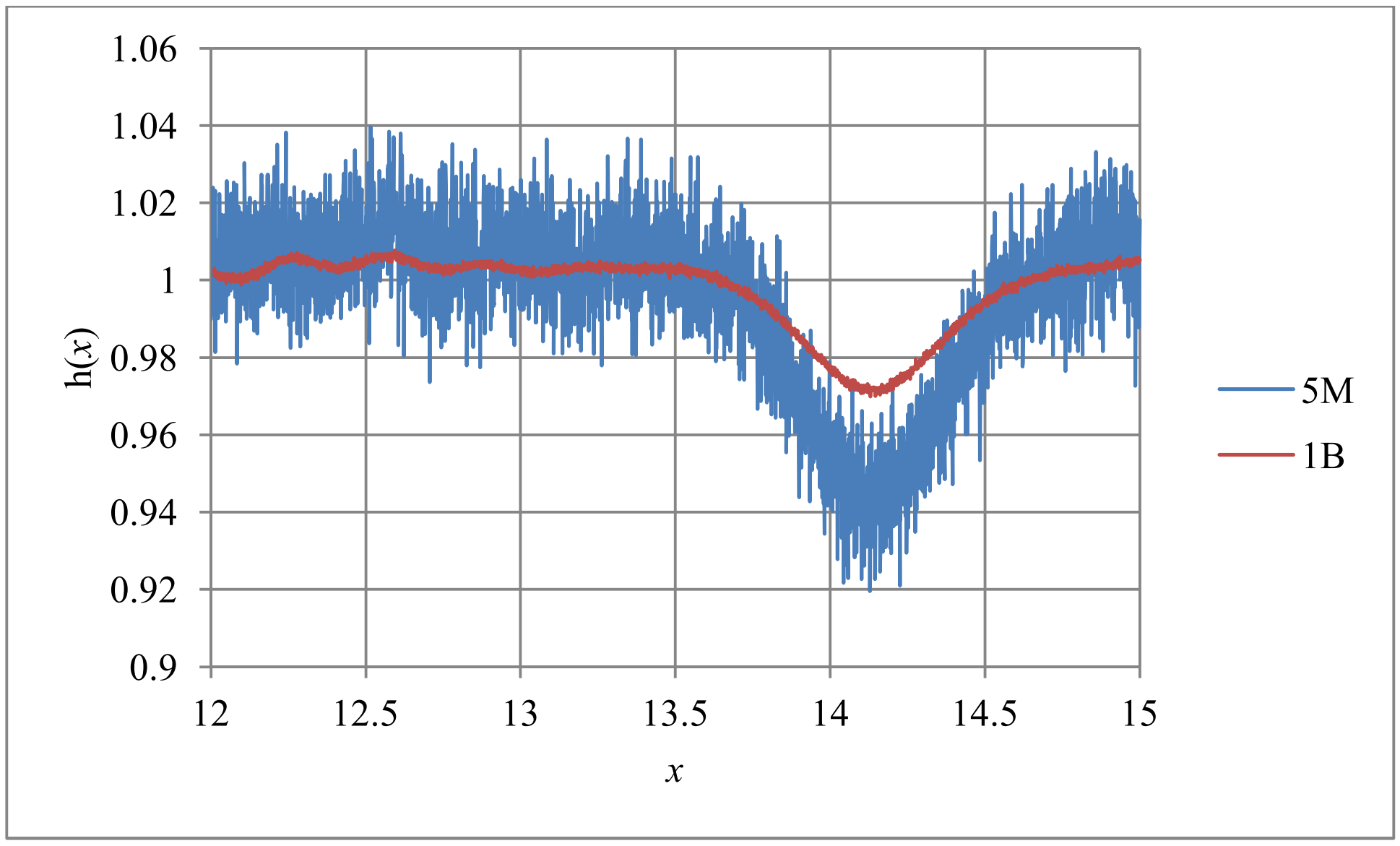

Fig. 1

The surrounding of the first Riemann zero =14.1347251… clearly shows the effect claimed and demonstrates that the amplitude of the repulsion diminishes with the increase of the amount of the Riemann Zeros used for statistics. The goal of this paper is to provide the exact "up to multiplicative constant" formula for this effect.

Fig. 2 represents lower portion of the spectra of the above two histograms[2].

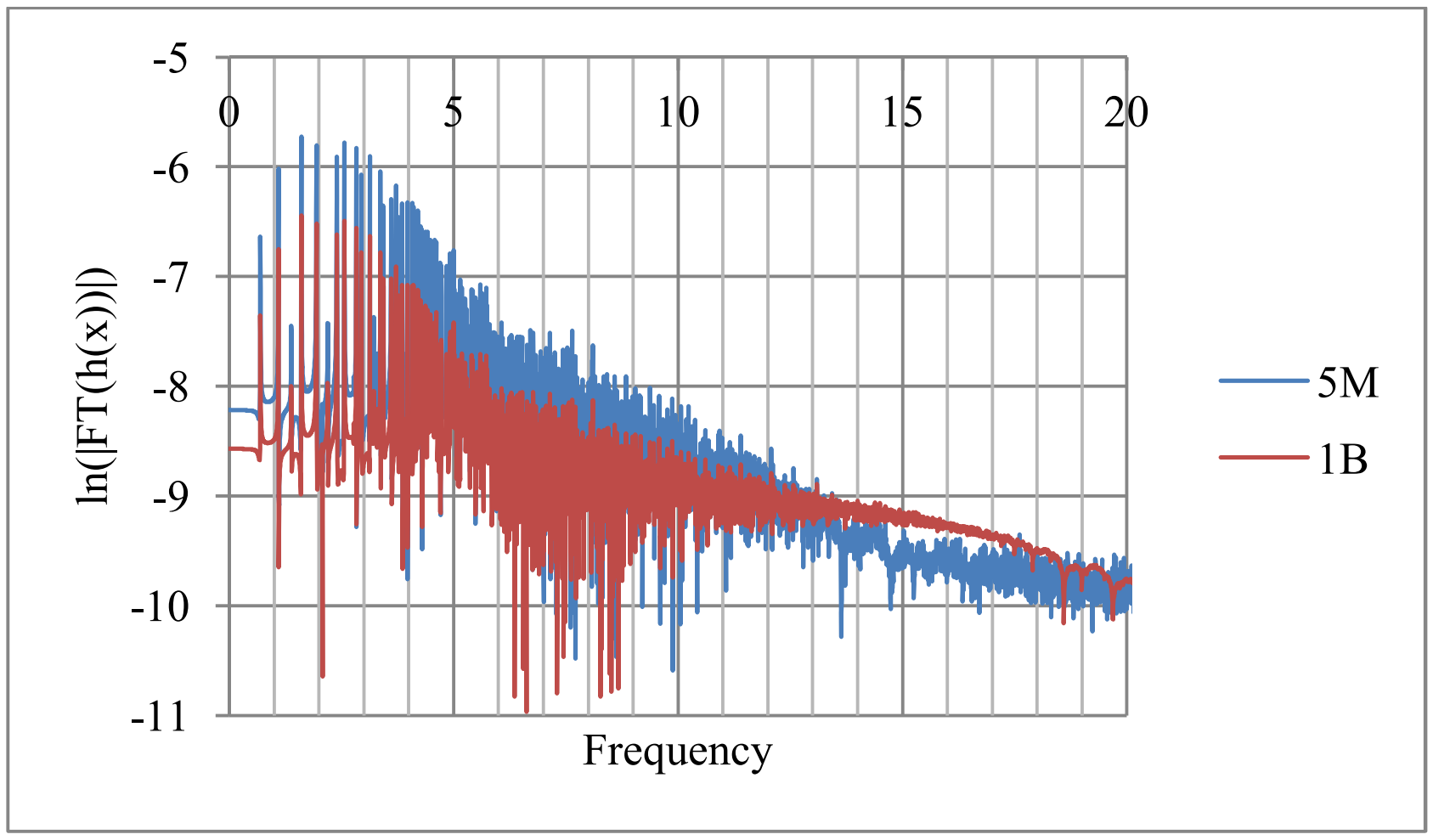

Fig. 2

---

[1] Possible only due to the great effort of David Platt published at http://www.lmfdb.org . The author owes him sincere appreciation.

[2] Following the example of Digital Library of Mathematical Functions http://dlmf.nist.gov let's use the designation $\ln(x)$ for natural (Napierian) logarithms.

In this paper number of graphs should be expected by the reader. One comment has to be made about the data displayed. Fig. 1 contains artificially scaled histograms in order to present them in the same scale. All Fourier Transform graphs contain the data scaled in such a way that $\mathrm{FT}(frequency = 0) = 1$. The value at zero frequency itself usually not displayed. All Fourier Transforms are made on non-symmetrized data, so evenness of the histograms is used only occasionally in the derivation. It means, among others, that the Imaginary part of the Fourier Transform is same relevant as its Real part.

These spectra (at least their low parts) clearly demonstrate presence of the (foreign to the original grid of spacing, which has been equal 0.001) frequencies multiple to $\ln p$ for all prime numbers $p$ (starting from $\ln 2 = 0.693 \ldots$). Fitting the multiples of logarithms of the lower primes it is easy to observe that the amplitude at $n\ln(p)$ is proportional to $p^{-n+1}$ for integer $n$. It means decomposition over

$$\sum_{n=1}^{\infty} \frac{e^{-in\ln(p)x}}{p^{n-1}} = \frac{p^{-ix}}{1 - p^{-1-ix}}$$

(or at least real part of it, because the histograms are even functions). This way terms

$$p \cdot \frac{1 - p\cos(x\ln p)}{p^2 - 2p\cos(x\ln p) + 1}$$

for all prime $p$ appear naturally. The actual form of this decomposition term indicates that we may be investigating some function on the line $\mathrm{Re}(z) = 1$.

The amplitudes of the terms (1) for various primes $p$ behave $\sim \frac{\ln^2(p)}{p}$, which may be seen from the envelope of the maxima on the Fig. 2.

Thus far we have obtain function

$$g_P(x) = \sum_{\text{prime } p \le P} \ln^2 p \cdot \frac{1 - p\cos(x\ln p)}{p^2 - 2p\cos(x\ln p) + 1}$$

These function appear in the histogram Fig. 1 with an amplitude depending only on the number of the Riemann Zeros [6] (RZs) used to create the histogram (and not depending on the prime $P$).

This dependency is tabulated in the Table 1 for first RZs between 5M and 10B (small portion of the David Platt's database, actually.) The increase of the correction amplitude with number of RZs is caused by the fact of increasing of total number of the points in the histogram (named original) with number of RZs. If the histograms are normalized to the same level (as done on Fig. 1), the correction amplitude decreases with increasing number of RZs.

| | | | | | | | | |
|---|---|---|---|---|---|---|---|---|
| Number of RZs | 5.0e+6 | 1.0e+7 | 1.5e+7 | 2.0e+7 | 2.5e+7 | 3.0e+7 | 3.5e+7 | 4.0e+7 |
| Amplitude | 133.241 | 252.916 | 368.366 | 481.233 | 592.250 | 701.830 | 810.250 | 917.692 |
| Number of RZs | 4.5e+7 | 5.0e+7 | 5.5e+7 | 6.0e+7 | 6.5e+7 | 7.0e+7 | 7.5e+7 | 8.0e+7 |
| Amplitude | 1024.291 | 1130.154 | 1235.359 | 1339.980 | 1444.063 | 1547.664 | 1650.814 | 1753.545 |
| Number of RZs | 8.5e+7 | 9.0e+7 | 9.5e+7 | 1.0e+8 | 2.0e+8 | 3.0e+8 | 4.0e+8 | 5.0e+e8 |
| Amplitude | 1855.895 | 1957.880 | 2059.527 | 2160.856 | 4139.075 | 6058.541 | 7941.685 | 9798.678 |
| Number of RZs | 1.0e+9 | 2.0e+9 | 3.0e+9 | 4.0e+9 | 5.0e+9 | 6.0e+9 | 7.0e+9 | 8.0e+9 |
| Amplitude | 18839.164 | 36271.986 | 53243.574 | 69928.589 | 86406.609 | 102723.273 | 118908.076 | 134981.636 |
| Number of RZs | 9.0e+9 | 1.0e+10 | | | | | | |
| Amplitude | 150959.162 | 166852.301 | | | | | | |

Table 1.
Amplitudes of the correction to the original histograms (for spacing 0.001),

Once the sum over primes in the definition of $g_P(x)$ is seemingly not convergent at $P \to \infty$, few prime values were chosen for $P$ to evaluate corrected histograms, namely $P$ = 4090441, 10001963, 603358771, and 1011890441. This choice corresponds to $\ln P$ =15.224…, 16.118…, 20.218…, and 20.735… correspondingly. The Fourier Transforms of $h_N(x) + A_N g_{4090441}(x)$ are presented on Fig. 3. This time real and imaginary parts of the Fourier Transform plotted together.

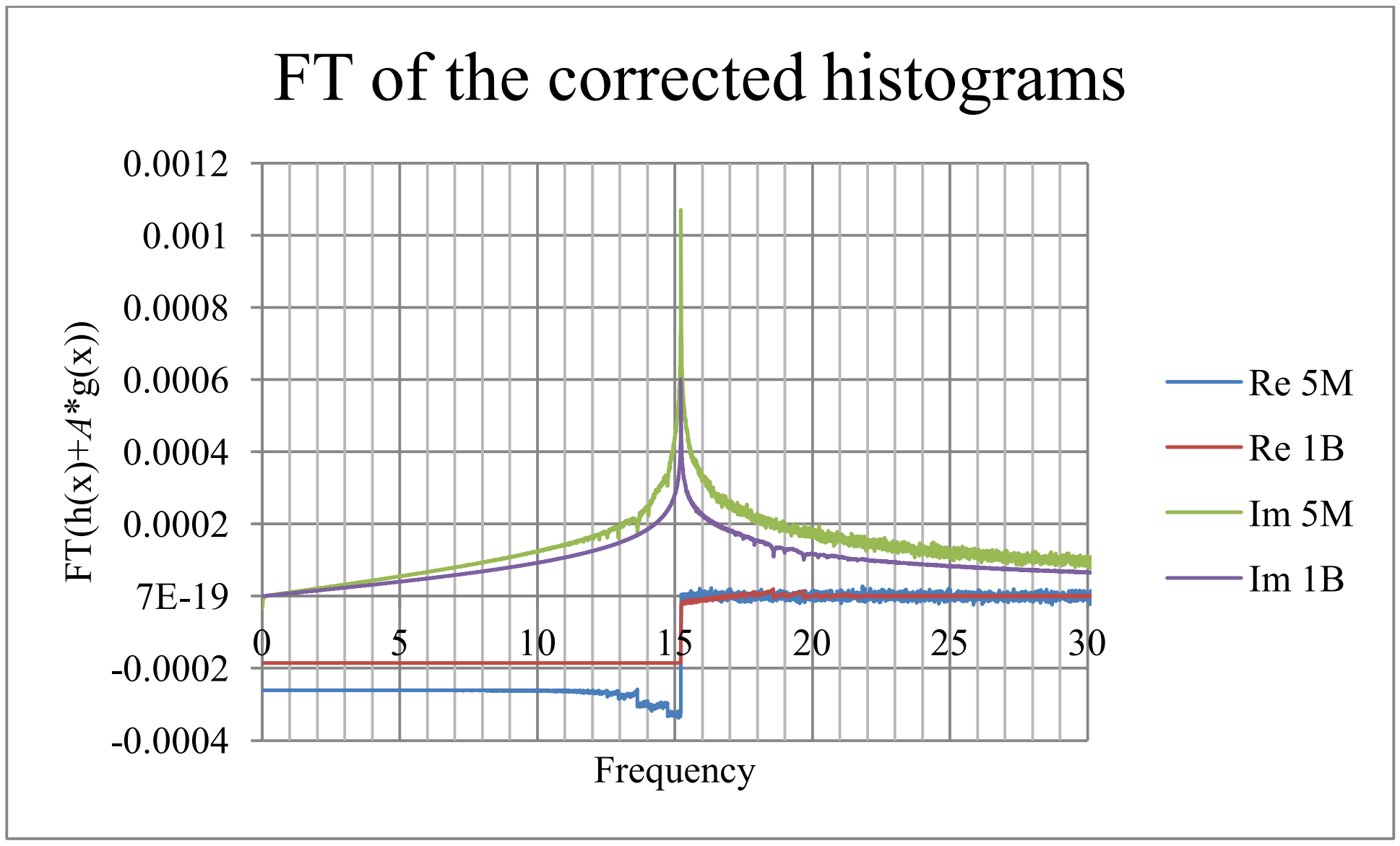

Fig. 3

The relative flatness of the real part till the frequency ~15.22… (and the corresponding log-type behavior of the imaginary part) tells that the erasure of the fixed frequencies has been achieved. Table 1 explains the absence of the exact formula (certainly not elementary) that still undeveloped.

The evident presence of the “noise” in the real (and imaginary) parts of the 5M curve might generate suspicions. Look at the Fig. 2, please. There are certain differences in the Fourier amplitudes, and they (amplitudes) may be sub-divided into two classes. The peaks we were eliminating so far are at fixed frequencies (multiples of the logarithms of the primes). However, there is some other "stuff" that does not fit perfectly with these set of fixed frequencies. The author thinks that it belongs to much bigger object [2], which he is not in a position to explain at this point. This "stuff" is moving with increasing number of the Riemann zeros used, − that is the reason why for 1B we do not see it on Fig. 3.

What we observe now is the jump of the real part near the limiting frequency $15.224\ldots = \ln 4090441$. This sudden jump may be easily corrected with the sinc-function. The amplitude of it is $A_N \pi \ln P$, that means that the definition of $g_P(x)$ should be corrected and now reads

$$\tilde{g}_P(x) = \ln P \frac{\sin(x \ln P)}{x} + \sum_{\text{prime } p \leq P} \ln^2 p \cdot \frac{1 - p\cos(x \ln p)}{p^2 - 2p\cos(x \ln p) + 1}$$

Additional term improves oscillations the sum over primes has for small $x$, but does not remove them completely. Careful examination of this behavior requires addition of the more rapidly decaying term, that has the same frequency $\ln P$. Not surprisingly enough this term leads to the very remarkable definition:

$$f'_P(x) = \ln P \cdot \frac{\sin(x \ln P)}{x} - \frac{1 - \cos(x \ln P)}{x^2} + \sum_{\text{prime } p \leq P} \ln^2 p \cdot \frac{1 - p\cos(x \ln p)}{p^2 - 2p\cos(x \ln p) + 1} \tag{1'}$$

## The Function

It is not occasional, of course, that the definition of the function we are interested in contains the prime, − it is a complete derivative:

$$f_P(x) = \frac{1 - \cos(x \ln P)}{x} - \sum_{\text{prime } p \leq P} \ln p \cdot \operatorname{atan} \frac{\sin(x \ln p)}{p - \cos(x \ln p)} \tag{1}$$

No claim of convergence of this definition as $P \to \infty$ is currently made, though the numerical evidence supports such thoughts.

Fig. 4 demonstrates the behavior of $f_P(x)$ at small values of $x$. It has been obtained with $P = 603358771$. The values of the function in the calculated region $[0{:}\,1000]$ belong to $[-2{:}\,2]$ and hopefully even to $[-\pi/2 : \pi/2]$.

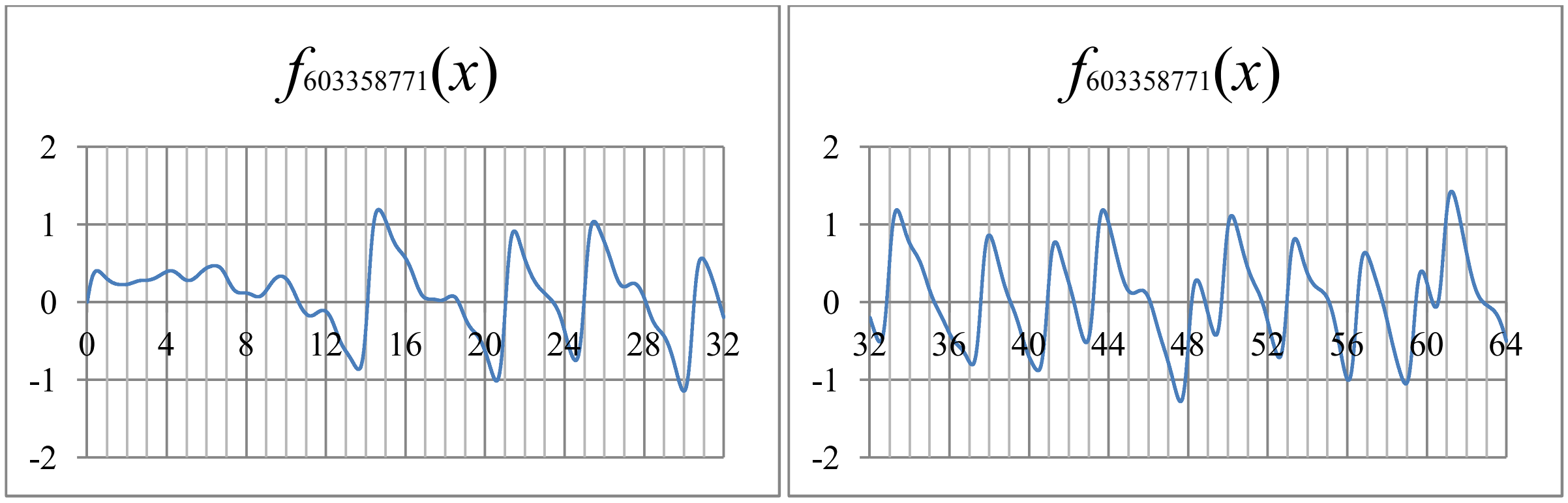

Fig. 4

To demonstrate the numerical convergence (better say, the lack of divergence) the difference of $f_P(x)$ for various $P$s is provided on Fig. 5a; analogous difference for $f'_P(x)$ provided on Fig. 5b.

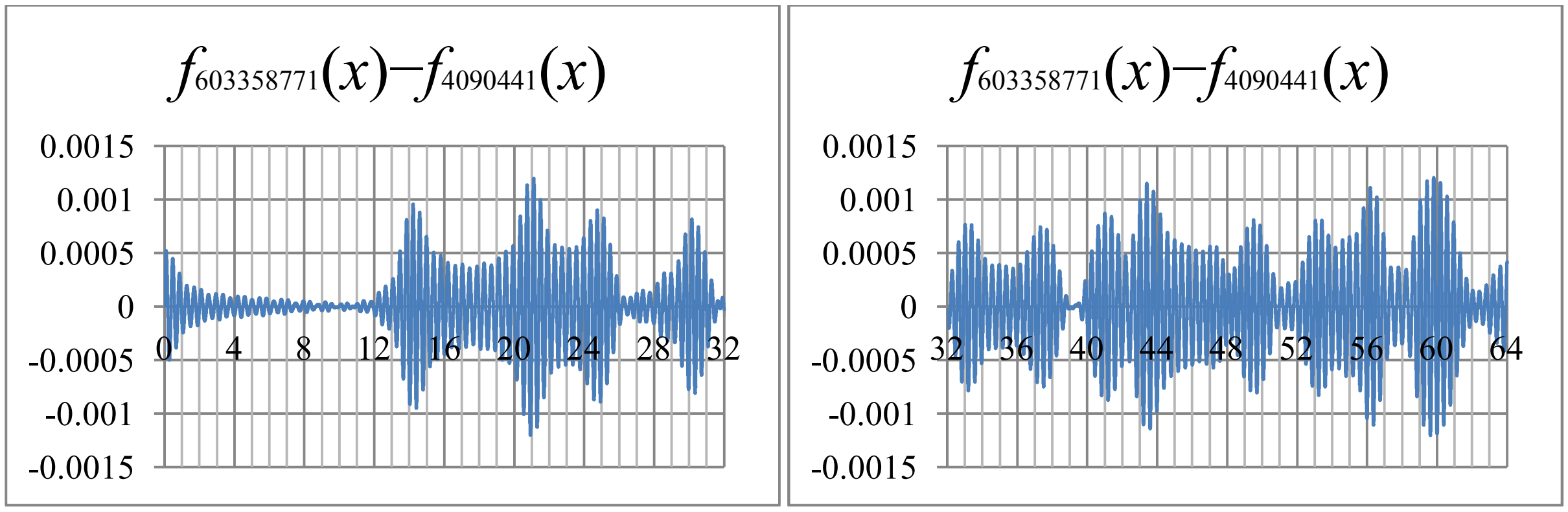

Fig. 5a

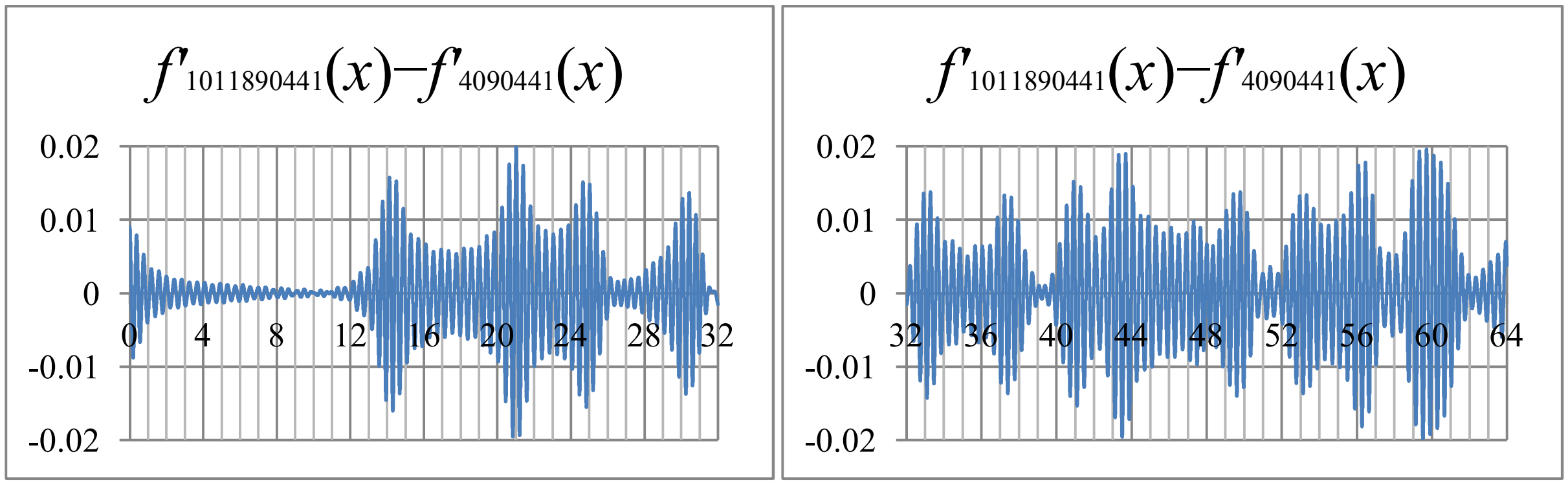

Fig 5b

Although Fig. 4 slightly resembles behavior of the Arg $\zeta(1/2+ix)$, it visibly does not coincide with it.

We started with exploration of the difference histograms of the RZs, and found that the part we are interested in is mostly described by the function (1′). Though the pictures provided so far relate only to $f_P(x)$, and not $f'_P(x)$. Fig. 6 covers this drawback presenting derivatives for $P=4090441$ (It does not matter what prime to use, because the graph resolutions do not allow us to see the differences. Unfortunately, Excel is not capable to plot more than 32K points.) The positions of the Riemann Zeros are marked on the *x*-axis.

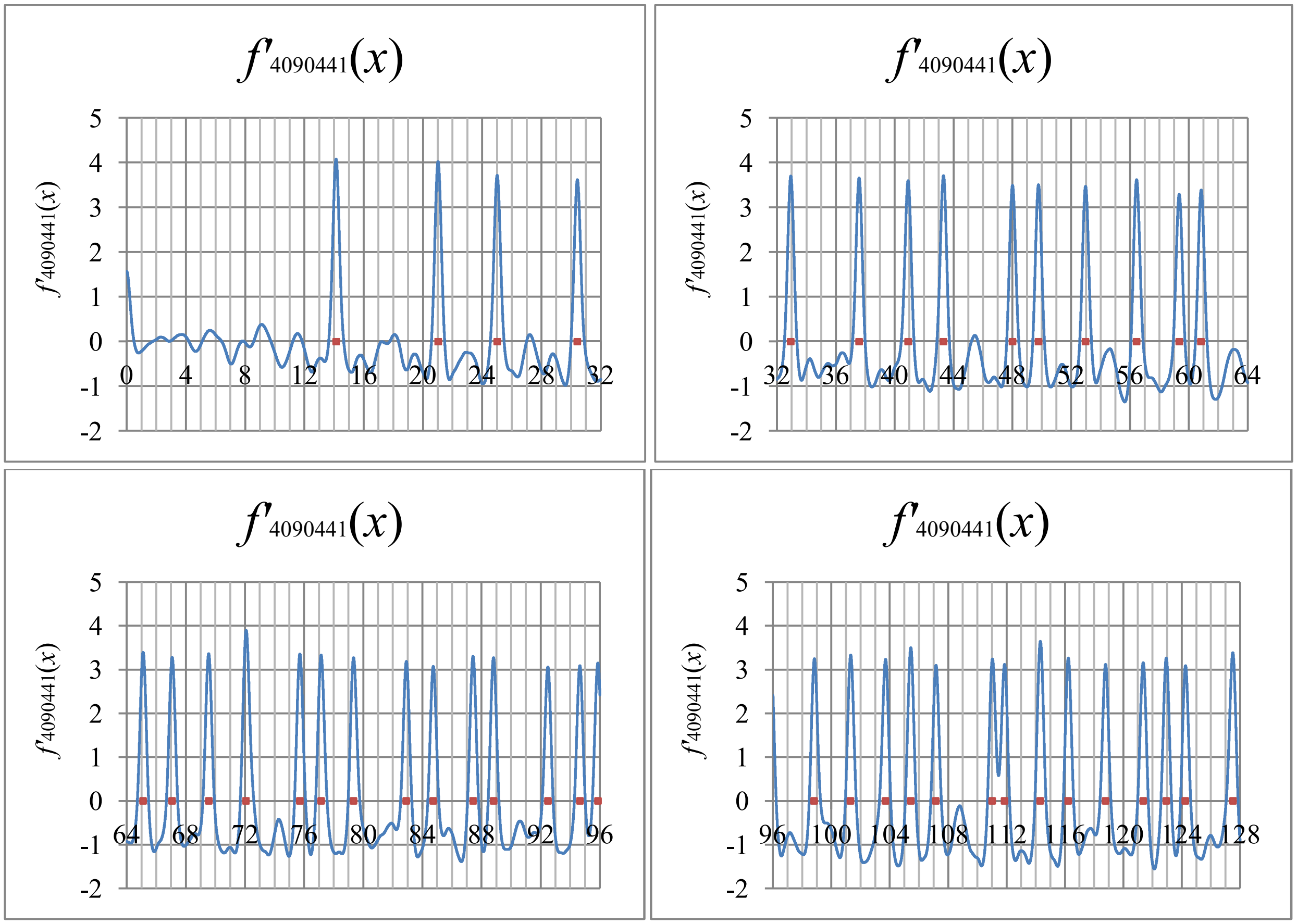

Fig. 6

The derivative for the calculated region $[0: 1000]$ always belongs to $[-2: 5]$. The positions of the RZs correspond to the values of the derivative $\geq 2$.

## Discussion

At this point reader should rebel: the referenced work [1] contains much more than just statistics of the RZs differences! Why are you limiting us with this part only, and do not even mention other *L*-functions, that do know about RZs too? Where is the Gaussian Unitarian Ensemble [5] references and consequences? What about more recent developments of the quantum chaos [3]?

First of all, the Manuscripts regarding eñe product, intensely referred in [1], are unavailable.

Second, the other *L*-functions have not been explored yet for the luck of the time and computer power. The major goal of this paper is to start discussion and probable thinking about "seemingly divergent" series.

Third, GUE [5] is, of course, present. It is in the remaining pieces of the histogram after correction with the function (1′) (with corresponding amplitude.) For better description, let's continue with the graphs. Figs. 7a−7b represent Fourier transforms of the corrected data for multiple number of RZs histograms (hope the evident designations do not confuse the reader.)

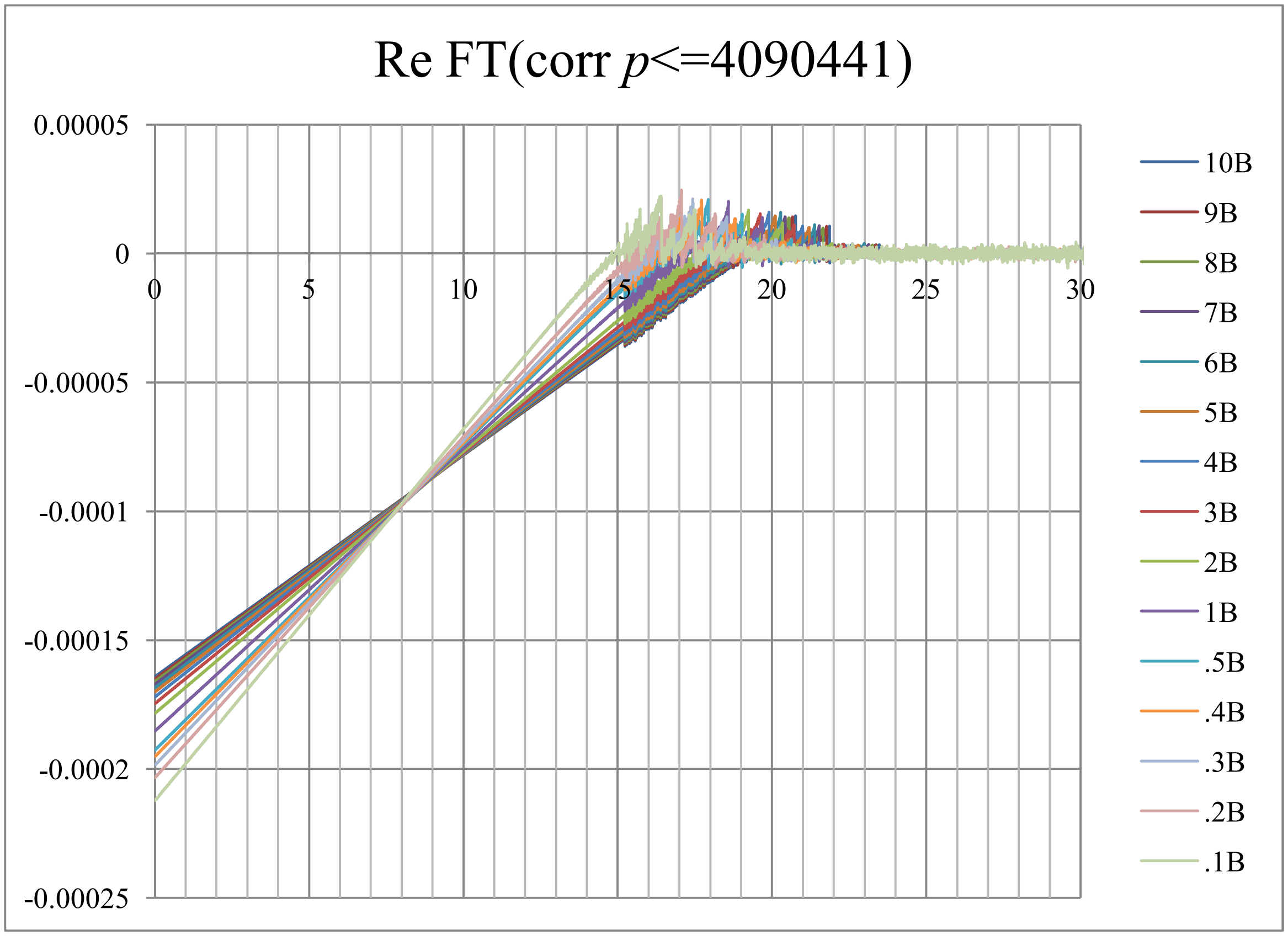

Fig. 7a

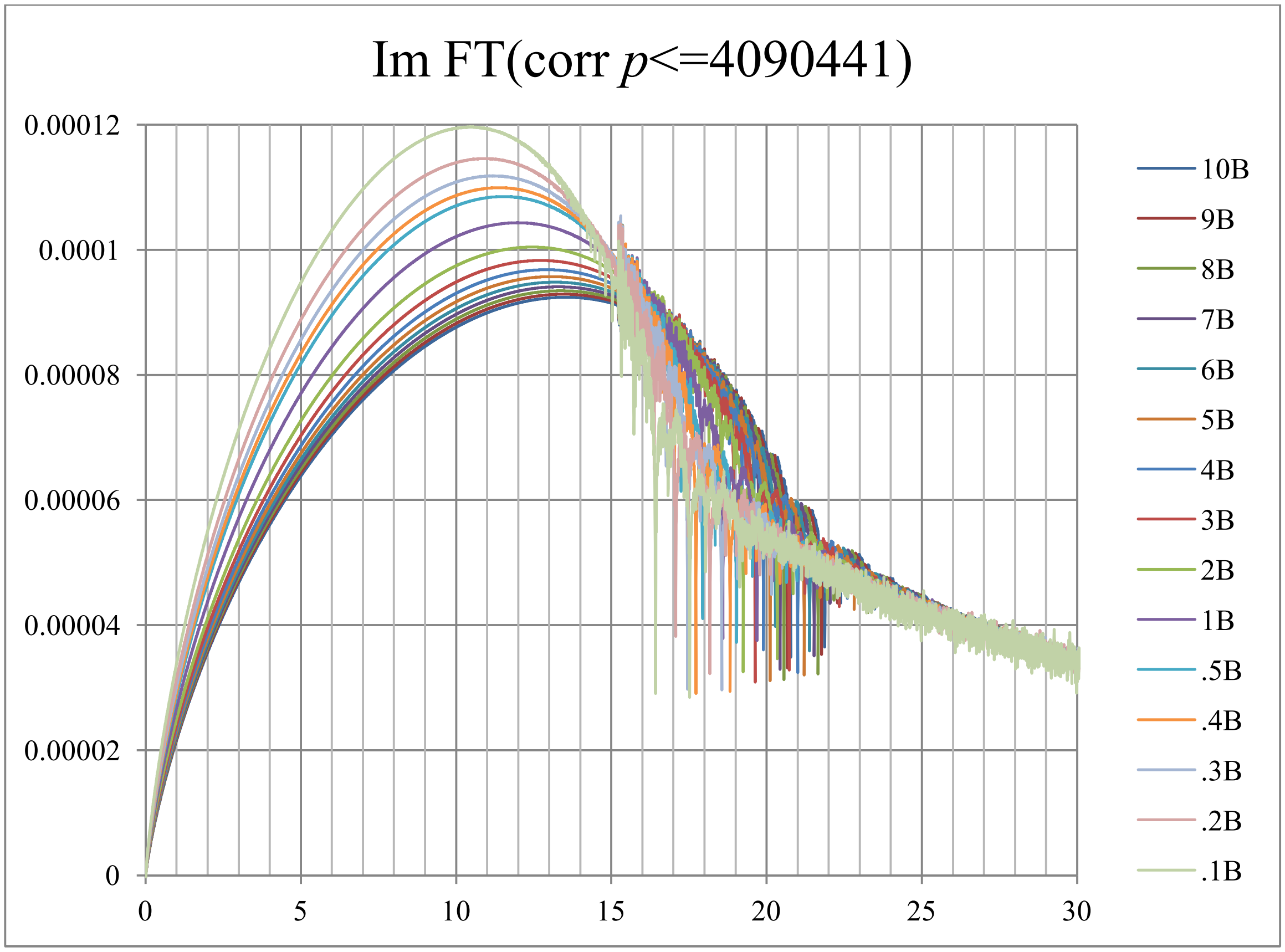

Fig. 7b

Plots of the remainders themselves are skipped in this paper because the excessiveness of the illustrative information.

The first term of (1) is just the derivative of $\mathrm{Cin}(x)$ in [5] and http://dlmf.nist.gov. The GUE behavior explained by the terms we have added to get rid of the oscillations of the $g_P(x)$.

The fact that interesting part of the histograms is described by the derivative (1′) is convenient in other way. Dealing with histograms we have always be aware of the decimation. At any point the consecutive points of the histogram may be summed together increasing the effective histogram spacing. It is very convenient that this "integration" leads to the boundary differences of The Function (1) itself, and might be evaluated as its derivative at some intermediate point (with much less restrictive requirements for the function itself) multiplied by a "bigger spacing". It means, that amplitudes tabulated in Table 1 may be considered as spacing proportional.

The straight lines in the graph Fig. 7a are evidently rotating, when the number of RZs involved in the histograms increasing. The center of rotation, may be at the frequency $2 \cdot \ln 73 = 8.5809 \ldots$.

What about quantum chaos? More presumably it includes comparison of (1) illustrated on Fig. 4 with the "moving stuff" that begins to appear at low numbers of RZs on Fig. 3. The reader should be aware that the functions suggested for comparison obtained from the different spaces (direct $x$ and frequencies), and if any similarity found, it indicates quantum replication phenomena.

Among the many questions remain open are the formulas like

$$\frac{\ln^2 P}{2} - \sum_{\text{prime } p \le P} \frac{\ln^2 p}{p-1} \xrightarrow[P\to\infty]{} \text{const} = 1.572\ldots$$

that confirms the value $f'(0)$[1], and many others. The last constant, if any, is probably different from $\pi/2$ and rather similar to the one found in [4]. The last work contains expressions similar (though evidently convergent) to the prime sum above.

What can we say about convergence of (1) as $P \to \infty$? First, since

$$\operatorname{atan}\frac{\sin(t)}{p-\cos(t)} = \frac{1}{2i}\ln\frac{p-e^{-it}}{p-e^{it}} = \sum_{n=1}^{\infty}\frac{\sin(nt)}{np^n}$$

and we may, if we wish, to ignore terms with $n \ge 2$, since sum over primes converges for them. It leads to

$$f_P(x) = \frac{1-\cos(x\cdot\ln P)}{x} - \sum_{\text{prime } p\le P}\sum_{n=1}^{\infty}\ln p\cdot\frac{\sin(nx\cdot\ln p)}{np^n}$$

Integrating we obtain

$$F_P(x) \coloneqq \int_0^x f_P(t)dt = \int_0^{x\cdot\ln P}\frac{1-\cos t}{t}dt - \sum_{n=1}^{\infty}\sum_{\text{prime } p\le P}\frac{1-\cos(nx\cdot\ln p)}{n^2p^n}$$

This expression may contain divergent constant, once we have to differentiate to obtain $f(x)$, though it is certainly less than $\ln(\ln P)$ [2]. The first term of $F_P(x)$ is $\operatorname{Cin}(x\cdot\ln P)$.

---

[1] The last formula is very similar to $\sum_{\text{prime } p\le x}\ln p/p = \ln x + O(1)$, published at http://dlmf.nist.gov/27.11.E10, but the presence of the square of logarithm, not the first power, changes the attitude.

[2] http://dlmf.nist.gov/27.11.E8.